\magnification=\magstephalf
\hfuzz=15pt
\baselineskip=14pt

\input xy
\xyoption{all}


\input amssym.def
\input amssym.tex

\def\build#1_#2^#3{\mathrel{\mathop{\kern 0pt#1}\limits_{#2}^{#3}}}

\def\Box{\square}
\def\geq{\geqslant}
\def\leq{\leqslant}

\def\ra{\rightarrow}
\def\lra{\longrightarrow}

\def\lla{\longleftarrow}

\def\egr{\build\equiv_{r}^{}}

\def\PP{\bf P}
\def\l{\ell}

\def\ip{{\frak P}}

\def\o{\omega}

\def\s{\sigma}
\def\d{\delta}

\def\a{\alpha}

\def\K{{\frak R}}
\def\J{{\frak A}}
\def\JJ{{\frak B}}

\def\Om{\Omega}
\def\Oo{{\cal O}}
\def\om{\omega}
\def\rk{\hbox{\rm rk}\, }
\def\hom{\hbox{\rm Hom}}
\def\sym{\hbox{\rm Sym}}
\def\codim{\hbox{\rm codim}\, }
\def\depth{\hbox{\rm depth}\, }
\def\he{\hbox{\rm ht}\, }
\def\ip{{\frak p}}
\def\ia{{\frak A}}
\def\lra{\longrightarrow}

\def\isoo{\build\cong_{1}^{}}

\def\isof{\build\cong_{5}^{}}

\def\cg{\lbrack\!\lbrack}
\def\cd{\rbrack\!\rbrack}
\def\ra{\rightarrow}
\def\lra{\longrightarrow}
\def\o{\omega}
\def\s{\sigma}

\def\d{\delta}
\def\eqr{\build\cong_{r}^{}}

\def\tc{{true codimension\ }}
\font\bigrm=cmb10 scaled 1100


\ \bigskip\bigskip\bigskip
\centerline{\bigrm HILBERT SERIES OF RESIDUAL INTERSECTIONS}\bigskip 

\centerline{Marc Chardin, David Eisenbud\footnote{$^{\diamondsuit}$}{partially
supported by 
the National Science Foundation.}, and Bernd Ulrich$^{\diamondsuit}$
\footnote{$^\heartsuit$}{ We are all grateful to MSRI for supporting the commutative algebra years in 2002-03 and 2012-13, where we worked on this paper.}
\footnote{}{AMS 2010 Subject Classification: Primary, 13D40, 13C40;  Secondary, 13H15, 13M06, 14C17, 14N15}}

\bigskip\bigskip
\noindent{\bf Abstract} We find explicit formulas for the Hilbert series of residual intersections of a scheme in terms of the Hilbert series of its conormal modules. In a previous paper we proved that such formulas should exist. We give applications to the dimension of secant varieties of surfaces and three-folds.
\bigskip

\noindent{\bf Introduction}\medskip

\noindent Let $M=\oplus_{i\in \bf Z} M_i$ be a finitely generated
graded module over the homogeneous coordinate ring $R$
of a projective variety over a field $k$. The {\it Hilbert series} (sometimes
called the Hilbert-Poincar\'e series) of $M$, which we will denote
$\cg M\cd$, is the Laurent series
$$
\cg M\cd = \sum (\dim M_i)t^i.
$$
If $Z\subset {\PP}^n:={\PP}^n_k$ is a scheme, then the Hilbert series of $Z$
is by definition the Hilbert series of the homogeneous coordinate
ring of $Z$. Of course this Hilbert series contains the data of the 
Hilbert polynomial of $Z$ as well.

Sometimes interesting geometric data (such as the dimension
of a secant variety) can be described in terms of {\it residual
intersections} in the sense of Artin and Nagata [2], and the purpose
of this paper is to compute the Hilbert series of such
schemes. Here is the 
definition: let $X\subset Y\subset {\PP}^n$ be closed subschemes 
of ${\PP}^n$,  let $R$ be the homogeneous coordinate ring of $Y$, and let
$I_X\subset R$ be the ideal of $X$ in $Y$. A
scheme $Z\subset Y$ is
called an $s$-residual intersection of $X$ in $Y$
if $Z$ is defined by an ideal of the form
$I_{Z}=(f_{1},\ldots ,f_{s}):_{R}I_{X}$, with $f_1,\dots,f_s$ homogeneous
elements in $I_X$, and
$Z$ is of 
codimension at least $s$ in $Y$.

We wish to derive formulas for the
Hilbert series
of $Z$ in terms of information about $X$
and the degrees of the polynomials $f_{i}$. 
In our previous paper [6] we showed that
this is sometimes possible in principle:
under certain hypotheses  the Hilbert series of $Z$ 
does not vary if we change the polynomials $f_i$, keeping
their degrees fixed.
In this paper we make this more precise by 
giving formulas---under somewhat stronger 
hypotheses---for the Hilbert series of $Z$ in terms of the
degrees of the $f_{i}$ and the Hilbert series of finitely many 
modules of the form $\o_R/I^m_{X}\o_R$, 
where $\o_R$ denotes the canonical module of $R$. 

For example, suppose that $Y=\PP^{n}$ and $X$ is locally a complete intersection (for instance, smooth).
If $f_{1}\dots, f_{s}$ are homogeneous elements of degree $d$ of $I=I_{X}$ such that $\K:=(f_{1}, \dots, f_{s}):I$ has codimension $\geq s$ then
the Hilbert series of the homogeneous coordinate ring of the scheme $Z$ defined by $\K$ differs from that
of a complete intersection defined by $s$ forms of degree $d$ by
$$
 \sum_{j=g}^{s}(-1)^{n+j}{s\choose j}t^{jd} \cg \o_{R}/I^{j-g+1}\o_{R} \cd (t^{-1}) +\hbox {a polynomial}.
$$
The polynomial
remainder term is present because we have made assumptions only on the scheme, and not on the
homogeneous coordinate ring. Here the expression $\cg R/I^{j-g+1}\cd(t^{-1})$ denotes the Laurent series obtained by writing 
$\cg R/I^{j-g+1}\cd$ as a rational function in $t$, substituting $t^{-1}$ for $t$, and rewriting
the result as a Laurent series. Up to a small shift in notation, this is the formula given in Remark 1.5b (where the
case of forms $f_{i}$ of different degrees is also treated.)

In applications, one sometimes needs to know only whether the 
$s$-residual intersection $Z$ actually has codimension exactly $s$
in $Y$;
for example, we will use such information in Section 3 to
say when the secant varieties of certain
(possibly singular) surfaces and smooth 3-folds have dimension less
than the expected dimension. For
this purpose it is enough to know just one coefficient of the 
Hilbert polynomial of
$Z$, that corresponding to the degree of
the codimension $s$ component of $Z$. More generally,
we show
how to use partial information about $X$ to compute
just the first $k$ coefficients of the Hilbert polynomial of $Z$.

Consider the case where
$Y$ is Gorenstein and $X$ is 
Cohen-Macaulay. Suppose, further, that
locally in 
 codimension $i< s$, the subscheme
 $X\subset Y$ can be defined by $i$ equations,
and that, for $j\leq s-g$, 
$$
{\rm depth}\ {\cal I}_{X}^{j}/{\cal I}_{X}^{j+1}
\geq \dim X-j.
$$
If $Z$ is any $s$-residual subscheme of $X$ in $Y$, then
the Hilbert polynomial of $Z$ may be written in terms of the
Hilbert polynomials of ${\cal I}_{X}^{j}$ for
$ j\leq s-g+1$ (the explicit formula is given in Theorem 1.9).
Moreover, if $X$ and $Y$ satisfy some of our
hypotheses only up to some codimension $r$, then the formula
 gives the first $r$ coefficients of the Hilbert polynomial of $Z$.

Our formulas are derived in Section 1, which is the technical heart of the paper. 
To prove them we need
to adapt the arguments of Ulrich [19]. The delicate point is the use in that
paper of the isomorphism $R\simeq \o_{R}$ that holds for 
a Gorenstein ring. Since our rings are Gorenstein only up to a certain
codimension $r$, we know only that $\o_R$ is a line bundle
locally in codimension $r$, and this does not suffice to determine its
Hilbert series. Thus the module $\o_R$ must be brought into play.
For other work along these lines, see Cumming [7]. 

In case where $Y$ is Gorenstein,
the sheaves ${\cal I}_X^j/{\cal I}_X^{j+1}$ 
themselves play the crucial role in our formulas. If
$X$ were locally a complete intersection scheme, then
${\cal I}_X/{\cal I}_X^{2}$ would be a vector bundle and ${\cal I}_X^j/{\cal I}_X^{j+1}$
would be its $j$-th symmetric power, so it is
reasonable to hope that for ``nice'' ideals $I$ the Hilbert series
of the first
few conormal modules
should determine the rest.
We prove a general theorem of this kind in Section 2, and carry out the reduction
in some particular cases. For instance,
if $s={\rm codim}(X)$ then
the degree of an $s$-residual intersection scheme $Z$ in $Y={\Bbb P}^{n}$
may be calculated immediately from
 B\'ezout's Theorem: $\deg Z = (\prod_i \deg f_i) - \deg X$. This was
extended to a formula in the case
 $s={\rm codim}_Y X+1$ by St\"uckrad [18]
and to the case 
 $s={\rm codim}_Y X+2$ by Huneke and Martin [16]. 
 Our formula gives an answer in general, and we work this out
 explicitly for the case 
$s={\rm codim}_Y X+3$. 

In Section 3 we apply our results to the study of secant loci. Our general theorems imply
conditions in terms of Chern classes of a smooth embedded three-fold for the degeneracy of
the secant locus in terms of Chern classes and in terms of certain Hilbert coefficients. We also
recover the analogous criteria for surfaces with mild singularities, a case already by Dale [7a] and others.
\medskip\medskip

\goodbreak{\bf 1. Formulas for the Hilbert series of residual
intersections} 
\medskip

To deal with rings $R$ that are not equidimensional, we define the
{\it \tc }of a prime ideal $\ip\subset R$ to be $\dim R - \dim R/\ip$.
We say that an ideal $I$ satisfies the condition $*G_s$ if, for every
prime $\ip$   in $V(I)$ of \tc $<s$, the minimal number of
generators of $I_\ip$ is at most $ \dim R_\ip$ (the usual condition
$G_s$ is the same but for codimension instead of \tc).

{\bf Definition 1.1.} Let $R$ be a graded ring and $M,N$ graded
$R$-modules. We say that $M$ and $N$ are {\it equivalent up to \tc
$r$}, and write $M\eqr N$, if there exist graded $R$-modules
$W_{1},\ldots ,W_{n}$ with $W_{1}=M$, $W_{n}=N$ and homogeneous maps
$W_{i}\ra W_{i+1}$ or $W_{i+1}\ra W_{i}$, for $1\leq i\leq
n-1$, which are isomorphisms locally in \tc $r$.
A homogeneous
map which is an isomorphism up to \tc $r$ will be denoted by
$\build\lra_{r}^{\sim}$.\medskip 

Saying that $M\eqr N$ is of course much stronger than saying that
$M$ and $N$ are isomorphic locally at each prime of \tc 
$\leq r$. For example, any two modules $M$ and $N$  that
represent line bundles on
a projective variety of dimension $r$ satisfy the latter condition,
but  $M\eqr N$ implies that
they represent isomorphic line bundles! The need to provide explicit
maps that are locally isomorphisms in some \tc
between modules that are not in fact
isomorphic is what makes the work in this section delicate.

If $R$ is a Noetherian standard graded algebra over a field and $M$ a finitely 
generated
graded $R$-module, we will
denote the Hilbert series of $M$ by $\cg M\cd$. If
$M\eqr N$ then $\cg M\cd
-\cg N\cd$, written as a rational function, has a pole
of order less than $\dim R - r$ at 1; we will write this 
as $\cg M\cd \egr \cg N\cd$ and say that these
series are {\it $r$-equivalent}. 
Thus if $r= \dim R-1=:d$ then 
$\cg M\cd \egr \cg N\cd$ means that
the Hilbert polynomials
 of $M$ and $N$ agree, 
and in general if $r<\dim R$ then $\cg M\cd \egr \cg N\cd$ means that
the Hilbert polynomials of $M$ and $N$, written in the form
$$
a_d{d+t\choose d}+a_{d-1}{d-1+t\choose t}+\cdots, 
$$
 have the same coefficients
of $a_s$ for $ s\geq \dim R - 1 -r. $

We will also extend
the notation $\egr$  to arbitrary series that are rational functions with no
pole outside 1, by the same requirement, as soon as $\dim R$ is clear
from the context.

 The substitution $t\mapsto t^{-1}$ is a well-defined automorphism of the ring ${\bf Z}[t,t^{-1},(1-t)^{-1}]$,
since $(1-t)^{-1}=-t(1-t)^{-1}\in {\bf Z}[t,t^{-1},(1-t)^{-1}]$. 

\medskip
{\bf Lemma 1.2.} Let $R$ be a positively graded Noetherian algebra over a field
$k$. If for each prime $\ip$ of dimension $\geq \dim R -r$,
the ring $R_\ip$ is Cohen-Macaulay of dimension $\dim R-\dim R/\ip$, then
$$
\cg \omega_R\cd (t) \egr (-1)^{\dim R} \cg R\cd (t^{-1}).
$$
\medskip
{\it Proof.} Dualize a free resolution of $R$ over a polynomial ring $S$, and note that all the homology
modules other than ${\rm Ext}_{S}^{\scriptstyle {\rm codim}_{S}(R)}(R, \omega_{S})$ are supported in codimension $>r$.
$\Box$
\medskip
We next adapt some results of [6] and [20] to our context.\medskip

Suppose that $R$ is a local Cohen-Macaulay ring
and I is an ideal of height $g$. In the rest of this section we will often use the condition that 
 ${\rm depth}\, R/I^{j}\geq \dim
R/I -j+1$ for $1\leq j\leq s-g$.  These conditions are  satisfied if
 $I$ satisfies
$G_s$ and if moreover, $I$ has the sliding depth
property  or, more  restrictively, is strongly Cohen-Macaulay (which
means that for every $i,~$ the
$i$-th Koszul homology $H_i$ of a generating set $~h_1, \ldots, h_n~$
of $I$ satisfies   ${\rm depth} ~ H_{i} \geq \dim R-n+i~$ or is
Cohen-Macaulay, respectively)  ([13, 3.3] and [15, 3.1]).   The latter
condition always holds if $I$ is a Cohen-Macaulay almost complete
intersection  or a Cohen-Macaulay deviation 2 ideal of a Gorenstein
ring [3].  It is also satisfied for any ideal in the linkage class
of a complete intersection [14, 1.11]. Standard examples include
perfect ideals of grade 2 ([1] and  [9]) and perfect Gorenstein ideals of
grade 3 [21].

\medskip 

\medskip
{\bf Lemma 1.3.} {\sl Let $R$ be a finitely generated positively graded
algebra over a factor ring of a local Gorenstein ring and write $\o =\o_{R}$. 
Let $I$ be a homogeneous ideal of height
$g$, let $f_{1},\ldots ,f_{s}$ be forms contained in $I$ of degrees
$d_{1},\ldots ,d_{s}$, write 
$\J_{i}:=(f_{1},\ldots ,f_{i})$, 
$\J :=\J_{s}$, and
$\K_{i}=\J_{i}:I$. Assume that ${\rm ht}\, \K_{i}\geq i$ for $1\leq
i\leq s$ and ${\rm ht}\, I+\K_{i}\geq i+1$ for $1\leq i\leq s-1$. Further
suppose that, locally off $V(\J )$, the elements $f_{1},\ldots ,f_{s}$ form
a weak regular sequence on $R$ and on $\o$, and that locally in \tc
$r$ in $R$ along $V(\J )$, the ring $R$ is Gorenstein and
 ${\rm depth}\, R/I^{j}\geq \dim
R/I -j+1$ for $1\leq j\leq s-g$. Then~:\smallskip

{\rm (a)} $(R/\K_{i-1})(-d_{i})
\build\lra_{r}^{\sim}
 \J_{i}/\J_{i-1} $
 via multiplication by $f_i$ for $1\leq i\leq
s$.\smallskip

{\rm (b)} $0\ra (\o I^{j}/\o \J_{i-1}I^{j-1})(-d_{i})\build\lra_{}^{\cdot
f_{i}} \o I^{j+1}/\o \J_{i-1}I^{j}\lra \o I^{j+1}/\o \J_{i}I^{j}\ra 0$ is
a complex that is exact locally in \tc $r$ for $1\leq i\leq s$
and $\min \{ 1,i-g\} \leq j\leq s-g$.\smallskip

{\rm (c)} ${\rm Ext}^{i}_{R}(R/\K_{i},\o )\eqr (\o
I^{i-g+1}/\o \J_{i}I^{i-g})(d_{1}+\cdots +d_{i})$ for $0\leq i\leq s$, if
locally in \tc $r$ in $R$ along $V(\J )$, ${\rm depth}\, R/I^{j}\geq \dim R/I -j+1$
for $1\leq j\leq s-g+1$.}\medskip

{\bf Proof.} Adjoining a variable to $R$ and to $I$ and localizing, 
we may suppose that
${\rm grade}\, I>0$. Furthermore, our assumptions imply that $I$ is
$G_{s}$ and satisflies the Artin-Nagata condition $AN_{s-1}$ locally in \tc $r$; see [20, 
2.9(a)]. Likewise, in the setting of (c), $AN_{s}$ holds
locally in \tc $r$. \par

Part (a) holds along $V(\J)$ by [20, 1.7(g)],
and off $V(\J)$ because $f_1,\dots,f_i$ form a weak
regular sequence.
Moreover, the sequence of (b) is
obviously a complex and it is exact in \tc $r$ by [20, 2.7(a)].\par

For the proof of (c), we induct on $i$. For $i=0$, our assertion is
clear since ${\rm grade}\, I>0$ and therefore $R/\K_{0}=R$. Assuming that
the assertion holds for $R_{i}=R/\K_{i}$ for some $i$, $0\leq i\leq s-1$,
we are going to prove our claim for $R_{i+1}=R/\K_{i+1}$. To this end we
may suppose $r\geq i+1$. 

We first wish to prove that 
$$
\leqno{(1)}\quad\quad\quad \quad
{\rm Ext}^{i+1}_{R}(R_{i}/(f_{i+1}R_i :IR_i),\o )\eqr 
[ I\ {\rm Ext}^{i}_{R}(R_{i},\o )/f_{i+1}{\rm Ext}^{i}_{R}(R_{i},\o ) ](d_{i+1}).
$$
Using the exact sequence $$
0\ra R_i/(0:_{R_i}f_{i+1}R_i)(-d_{i+1})\buildrel{.f_{i+1}}\over{\lra} R_i\lra R_i/f_{i+1}R_i \ra 0
$$
we obtain a long exact sequence
$$
\cdots {\rm Ext}^{i}_{R}(R_{i},\o ) \buildrel{.f_{i+1}}\over{\lra} {\rm Ext}^{i}_{R}(R_i/(0:_{R_i}f_{i+1}R_i),\o )(d_{i+1})
\ra {\rm Ext}^{i+1}_{R}(R_{i}/f_{i+1}R_i,\o )\ra {\rm Ext}^{i+1}_{R}(R_{i},\o )\cdots
$$
Since  the support in $R$ of 
 $ (0:_{R_i}f_{i+1})$ has \tc $\geq r+1>i$
by part (a), we have 
${\rm Ext}^{i}_{R}(R_i/(0:_{R_i}f_{i+1}R_i),\o )\build\lra_{r}^{\sim}
{\rm Ext}^{i}_{R}(R_i,\o )$ via the natural map.
Furthermore ${\rm Ext}^{i+1}_{R}(R_{i},\o )\eqr 0$; as locally up to \tc $r$ on 
$V(\J )$, $R_i$ is Cohen-Macaulay of \tc $i$ by [20, 1.7(a)], whereas locally up to \tc $r$ off 
$V(\J )$, $R_i$ is defined by the weak regular sequence $f_1,\ldots, f_i$ and hence has projective dimension at most $i$. 
Therefore
$$
\left[ {\rm Ext}^{i}_{R}(R_i,\o )/f_{i+1}{\rm Ext}^{i}_{R}(R_i,\o )\right] (d_{i+1})\eqr 
E:={\rm Ext}^{i+1}_{R}(R_{i}/f_{i+1}R_i,\o ).
$$
Hence, to prove (1),  it suffices to show that 
$$
I\ {\rm Ext}^{i+1}_{R}(R_{i}/f_{i+1}R_i,\o )\eqr {\rm Ext}^{i+1}_{R}(R_{i}/(f_{i+1}R_i :IR_i),\o ).
$$
The map $R_i/f_{i+1}R_i\ra R_i/(f_{i+1}R_i:_{R_i}IR_i)$ induces a map
$$
\phi: \ {\rm Ext}^{i+1}_{R}(R_{i}/(f_{i+1}R_i :IR_i),\o )\ra {\rm Ext}^{i+1}_{R}(R_{i}/f_{i+1}R_i,\o ).
$$
We prove that locally in \tc $r$ in $R$, the map $\phi$ is injective,
and its image coincides with $I\;E$, which gives
$$
{\rm im}\ \phi \build\lra_{r}^{\sim} {\rm im\ } \phi+I\;E \build\lla_{r}^{\sim} I\;E.
$$
This is trivial locally off $V(\J )$ because on this locus $I=R$ and 
$f_{i+1}R_i :IR_i=f_{i+1}R_i$. 
Therefore, we may localize to assume that $R$ is a local ring of dimension at most $r$ and $\J \not= R$.
Of course we may suppose $R_i\not= 0$. In this case $R$ is Gorenstein, $R_i$ is Cohen-Macaulay 
of codimension $i$, and $f_{i+1}$ is a non zerodivisor on $R_i$
by [20, 1.7(f)]. Let $S=R_i/f_{i+1}R_i$. 

The natural equivalence of functors ${\rm Ext}^{i+1}_R(\hbox{---}, \omega)
\simeq {\rm Hom}_S (\hbox{---}, \om_S)$ together with the exact sequence
$$
0\ra 0:_S IS\lra S \lra S/0:_S IS\lra 0
$$
yield a commutative diagram with an exact row
$$
\xymatrix{
&{\rm Ext}^{i+1}_R(R_i/(f_{i+1}R_i :IR_i ), \omega)\ar^(.75){\phi}[r]
\ar^{\simeq}[d]&E\ar^{\simeq}[d]&&\\
0\ar[r]& {\rm Hom}_S (S/0:_S IS, \om_S)\ar^(.65){\psi}[r]&\om_S \ar[r]&
 {\rm Hom}_S (0:_S IS, \om_S)\ar[r]&0\\}
$$
The last map is surjective because $S/(0:_SIS)=R_{i+1}$
by [20,1.7(f)] and $R_{i+1}$ is a maximal Cohen-Macaulay
$S$-module. Now $\phi$ is injective and the desired equality
${\rm im}\ \phi = I\;E$ follows once we have shown that 
${\rm im\ }\psi = I\om_S$. For this it suffices to prove
$$
\leqno{(2)}\quad\quad\quad\quad\quad\quad\quad\quad\quad\quad\quad\quad
 {\rm coker}\;\psi \simeq \om_S/I\om_S;
$$
for then $I\om_S\subset {\rm im}\;\psi$ and we have the natural
epimorphism of isomorphic modules
$\om_S/I\om_S \to {\rm coker}\; \psi$,
which is necessarily an isomorphism.
We first argue that $\om_S/I\om_S$ is a maximal Cohen-Macaulay
$S$-module. Indeed, [20, 2.7(c)] gives 
$\K_i\cap I^{i-g+2} = \J_i I^{i-g+1}$, 
which implies 
$$
\leqno{(3)} \quad\quad\quad\quad\quad\quad\quad\quad\quad
\J_iI^{i-g}\cap I^{i-g+2} = \J_i I^{i-g+1}.
$$
Hence by our induction hypothesis,
$$
I\om_S\simeq I^{i-g+2}/(\J_iI^{i-g}\cap I^{i-g+2}+f_{i+1}I^{i-g+1})
=I^{i-g+2}/J_{i+1}I^{i-g+1}.
$$
But the latter is indeed a maximal Cohen-Macaulay $S$-module
according [20, 2.7(b)].
Thus
$$
\om_S/I\om_S\simeq 
{\rm Hom}_S ({\rm Hom}_S (\om_S /I\om_S,\om_S), \om_S)
\simeq
{\rm Hom}_S ({\rm Hom}_S (S /IS, S), \om_S)
\simeq
{\rm Hom}_S (0:_S IS, \om_S).
$$
This completes the proof of (2), and hence of (1).

 Now 
$$
I{\rm Ext}^{i}_{R}(R_{i},\o )/f_{i+1}{\rm Ext}^{i}_{R}(R_{i},\o )\eqr (\o I^{i-g+2}/(\o \J_{i}I^{i-g}\cap \o
I^{i-g+2}+\o f_{i+1}I^{i-g+1}))(d_{1}+\cdots +d_{i})
$$ 
by our induction hypothesis, and using (3) one sees that

$$
\leqno{(4)} \quad\quad\quad
I{\rm Ext}^{i}_{R}(R_{i},\o )/f_{i+1}{\rm Ext}^{i}_{R}(R_{i},\o )\eqr 
(\o I^{i-g+2}/\o \J_{i+1}I^{i-g+1})(d_{1}+\cdots +d_{i}).
$$ 
On the other hand, $R_{i+1}\build\lra_{r}^{\sim} R_i /(f_{i+1}R_i :IR_i )$ according to [20, 1.7(f)], and hence
$$
\leqno{(5)} \quad\quad\quad
{\rm Ext}^{i+1}_{R}(R_{i+1},\o )\eqr 
{\rm Ext}^{i+1}_{R}(R_i /(f_{i+1}R_i :IR_i),\o ).
$$ 

Now combining (5), (1), (4) concludes the proof of part (c).
\hfill$\Box$
\bigskip

We write
$
\s_{m}(t_{1}, \dots, t_{s}) 
$
for the $m$-th elementary symmetric function.

\medskip
{\bf Theorem 1.4.} Let $R$ be a standard graded Noetherian algebra over a field.
Write
$n=\dim R$ and $\o =\o_{R}$, and let $I$ be a homogeneous ideal of height $g$ satisfying
$*G_{s}$. Let $f_{1},\ldots ,f_{s}$ be forms contained in $I$ of degrees
$d_{1},\ldots ,d_{s}$, write $\Delta_{s}:=\prod_{i=1}^{s}(1-t^{d_i})$,
$\J =(f_{1},\ldots ,f_{s})$, $\K =\J :I$, and  
assume that ${\rm ht}\, \K \geq s$.
For each prime of \tc $\leq r$ suppose: 

\item{$\bullet$} If
 $\ip \notin V(\J )$, then the elements $f_{1},\ldots ,f_{s}$ form
a weak regular sequence on $R_\ip$ and on $\o_\ip$.

\item{$\bullet$} If $\ip \in V(\J)$, then
the ring $R_\ip$ is Gorenstein of dimension  equal to the
\tc of $\ip$
and
 ${\rm depth}\, R_\ip /I^{j}_\ip\geq \dim
R_\ip/I_\ip -j+1$ for $1\leq j\leq s-g$. 

\smallskip
{\rm (a)} 
$$
\cg R/\J \cd (t)\egr \Delta_{s}\cg R\cd (t)-(-1)^{n-g}
\sum_{j=1}^{s-g}(-1)^{j}\s_{g+j}(t^{d_1},\ldots ,t^{d_s}) \cg \o /I^{j}\o \cd (t^{-1}). 
$$ 

{\rm (b)} If furthermore, locally in \tc
$r$ in $R$ along $V(\J )$,
 ${\rm depth}\, R/I^{s-g+1}\geq \dim
R/I -s+g$, then
$$
\cg R/\K \cd (t)\egr \Delta_{s}\cg R\cd (t)-(-1)^{n-g}
\sum_{j=1}^{s-g+1}(-1)^{j-1}\s_{g+j-1}(t^{d_1},\ldots ,t^{d_s}) \cg \o /I^{j}\o \cd (t^{-1}). 
$$ 
\medskip

{\bf Proof.} For $0\leq i\leq s$, write $\J_{i}=(f_{1},\ldots ,f_{i})$,
$\K_{i}=\J_{i}:I$. Applying repeatedly a general position argument (see [6, 2.5]), we may assume that
${\rm ht}\, \K_{i}\geq i$ and ${\rm ht}\, I+\K_{i}\geq i+1$ for $0\leq
i\leq s-1$.

We first notice that for $0\leq i\leq s$ and $i-g\leq j\leq s-g$, 
$$
\cg \o \J_{i}I^{j}\cd \egr \sum_{\l =1}^{i}(-1)^{\l +1}\s_{\l
}(t^{d_{1}},\ldots ,t^{d_{i}})\cg \o I^{j-\l +1}\cd ,
$$
which can be easily deduced from Lemma 1.3(b) using induction on $i$.

Now let $0\leq i\leq s-1$ for (a), or $0\leq i\leq s$ for (b),
respectively. Then by Lemma 1.3(c),
$$
\leqno{(1)}\quad\quad\quad \eqalign{\cg {\rm Ext}^{i}_{R}(R/\K_{i},\o )\cd (t) &\egr t^{-(d_{1}+\cdots
+d_{i})}(\cg \o I^{i-g+1}\cd -\cg \o \J_{i}I^{i-g}\cd )(t)\cr
   &\egr t^{-(d_{1}+\cdots +d_{i})} \sum_{\l
=0}^{i}(-1)^{\l}\s_{\l}(t^{d_{1}},\ldots ,t^{d_{i}})\cg \o I^{i-g-\l
+1}\cd (t).\cr}
$$

Locally in \tc $r$ in $R$, $R/\K_{i}$ is either zero or
Cohen-Macaulay of codimension $i$ by [20, 2.9 and 1.7(a)]. Therefore
$$
\cg {\rm Ext}^{i}_{R}(R/\K_{i},\o )\cd (t)\egr (-1)^{n-i}\cg R/\K_{i}\cd (t^{-1}),
$$ 
as can be easily seen by dualizing a homogeneous finite free resolution
of $R/\K_{i}$ over a polynomial ring and using Lemma 1.2. Now by Lemma
1.3(c) and (1), 
$$
\eqalign{
\cg R/\K_{i}\cd (t)&\egr (-1)^{n-i}\cg {\rm Ext}^{i}_{R}(R/\K_{i},\o )\cd (t^{-1})\cr
      &\egr (-1)^{n-i}t^{d_{1}+\cdots +d_{i}}\sum_{\l
=0}^{i}(-1)^{\l}\s_{\l}(t^{-d_{1}}, \ldots ,t^{-d_{i}})\cg \o I^{i-g-\l +1}\cd
(t^{-1})\cr
      &=(-1)^{n-i}\sum_{\l =0}^{i}(-1)^{\l}\s_{i-\l}(t^{d_{1}},\ldots
,t^{d_{i}})\cg \o I^{i-g-\l +1}\cd (t^{-1}).\cr}
$$
Therefore
$$
\leqno{(2)}\quad\quad\quad \cg R/\K_{i}\cd (t)\egr (-1)^{n-g}
\displaystyle\sum_{j=-g+1}^{i-g+1}(-1)^{j+1}
\s_{g+j-1}(t^{d_{1}},\ldots , t^{d_{i}})\cg \o I^{j}\cd (t^{-1}).
$$

Now part (b) follows since
$$
\eqalign{
\Delta_{i} \cg R\cd (t)\egr \Delta_{i}(-1)^{n}\cg \o \cd (t^{-1})
&=(-1)^{n-i}\sum_{\l =0}^{i}(-1)^{\l}\s_{i-\l}(t^{d_{1}},\ldots
,t^{d_{i}})\cg \o \cd (t^{-1})\cr
&=(-1)^{n-g}\sum_{j=-g}^{i-g}(-1)^{j}\s_{g+j}(t^{d_{1}},\ldots
,t^{d_{i}})\cg \o \cd (t^{-1}).\cr}
$$

To see (a) notice that by Lemma 1.3(a), 
$$
\cg R/\J_{i}\cd (t)=\cg R/\J_{i-1}\cd (t)-t^{d_{i}}\cg R/\K_{i-1}\cd (t)
$$
for $1\leq i\leq s$. Now by induction on $i$ using (2), one shows that
$$
\cg R/\J_{i}\cd (t)\egr (-1)^{n-g}\sum_{j=-g}^{i-g}(-1)^{j}\s_{g+j}
(t^{d_{1}},\ldots ,t^{d_{i}})\cg \o I^{j}\cd (t^{-1}),
$$
from which (a) can be easily deduced.\hfill$\Box$\medskip

{\bf Remark 1.5.} If in Theorem 1.4, $R$ is a polynomial ring in $n$
variables, then the formulas of that theorem take the following form:
$$
\leqno{\hbox{(a)}}\quad\quad\quad\quad \cg R/\J \cd (t)\egr \Delta_{s}\cg R\cd (t)-(-t)^{-n}\displaystyle
\sum_{j=1}^{s-g}(-1)^{g+j}\s_{g+j}(t) \cg R/I^{j} \cd (t^{-1}); 
$$ 
$$
\leqno{\hbox{(b)}}\quad\quad\quad\quad
\cg R/\K \cd (t)\egr \Delta_{s}\cg R\cd (t)-(-t)^{-n}\displaystyle
\sum_{j=1}^{s-g+1}(-1)^{g+j-1}\s_{g+j-1}(t) \cg R/I^{j} \cd (t^{-1}). 
$$

{\bf Lemma 1.6.} Write $\Delta_{s}(t)=(1-t)^{s}\sum_{k\geq
0}c_{k}(d_{1},\ldots ,d_{s})(1-t)^{k}$. Then
$$
c_{k}(d_{1},\ldots ,d_{s})=(-1)^{k}\sum_{{i_{1}\geq 1,\ldots ,i_{s}\geq
1}\atop{\scriptscriptstyle{i_{1}+\cdots
+i_{s}=k+s}}}\prod_{j=1}^{s}{{d_{j}}\choose{i_{j}}}.
$$

{\bf Proof.} Write $P_{j}(t)=\sum_{\l =0}^{d_{j}-1}t^{\l}$ and notice
that $\Delta_{s}(t)=(1-t)^{s}\prod_{j=1}^{s}P_{j}(t)$.

Now $P_{j}^{(m)}(1)=\sum_{\l =0}^{d_{j}-1}m!{{\l}\choose{m}}=m!
{{d_{j}}\choose{m+1}}$. Hence
$$
\eqalign{\Bigl( \prod_{j=1}^{s}P_{j}\Bigl) ^{(k)}(1)&=
\sum_{{m_{1}\geq 0,\ldots ,m_{s}\geq 0}\atop{\scriptscriptstyle{m_{1}+\cdots
+m_{s}=k}}}{{k!}\over{m_{1}!\cdots
m_{s}!}}\prod_{j=1}^{s}P_{j}^{(m_{j})}(1)\cr
       &=k!\sum_{{m_{1}\geq 0,\ldots ,m_{s}\geq 0} 
\atop{\scriptscriptstyle{m_{1}+\cdots +m_{s}=k}}}\prod_{j=1}^{s}
{{d_{j}}\choose{m_{j}+1}}.\cr}
$$
This yields our formula since $c_{k}(d_{1},\ldots
,d_{s})={{(-1)^{k}}\over{k!}}\left( \prod_{j=1}^{s}P_{j}\right)
^{(k)}(1)$.\hfill$\Box$\medskip

{\bf Lemma 1.7.} Let $P$ be a numerical polynomial written in the form
$P(t)=\sum_{i=0}^{m}(-1)^{i}e_{i}{{t+m-i}\choose{m-i}}$. For an integer
$d$ define the polynomial $Q(t)=P(-t+d)$, and write
$Q(t)=\sum_{i=0}^{m}(-1)^{i}h_{i}{{t+m-i}\choose{m-i}}$. Then
$$
h_{i}=(-1)^{m}\sum_{k=0}^{i}(-1)^{k}{{d+m+1-k}\choose{i-k}}e_{k}.
$$

{\bf Proof.} We first notice that for integers $r$ and $n\geq 0$, one
has the following identities of numerical polynomials~:
$$
\leqno{\hbox{(3)}}\quad\quad\quad\quad\displaystyle
{{-t+n}\choose{n}}=(-1)^{n}{{t-1}\choose{n}},
$$
$$
\leqno{\hbox{(4)}}\quad\quad\quad\quad\displaystyle
{{t+r+n}\choose{n}}=\sum_{\l =0}^{n}{{r-1+\l
}\choose{\l}}{{t+n-\l}\choose{n-\l}},
$$
where the first is obvious and the second can be easily proved by
induction on $n$.

Now 
$$
\eqalign{
Q(t)&=\sum_{k=0}^{m}(-1)^{k}e_{k}{{-t+d+m-k}\choose{m-k}}\cr
    &=(-1)^{m}\sum_{k=0}^{m}e_{k}{{t-d-1}\choose{m-k}}\quad {\rm by}\ (3)\cr
    &=(-1)^{m}\sum_{k=0}^{m}e_{k}{{t+(-d-1-m+k)+(m-k)}\choose{m-k}}\cr
    &=(-1)^{m}\sum_{k=0}^{m}e_{k}\sum_{\l
=0}^{m-k}{{-d-1-m+k-1+\l}\choose{\l}}{{t+m-k-\l}\choose{m-k-\l}}\quad {\rm by}\ (4)\cr
    &=(-1)^{m}\sum_{i=0}^{m}\left(\sum_{k=0}^{i}{{-d-m-2+i}\choose{i-k}}
e_{k}\right) {{t+m-i}\choose{m-i}},\cr}
$$
where
${{-d-m-2+i}\choose{i-k}}=(-1)^{i+k}{{d+m+1-k}\choose{i-k}}$
by (3).\hfill$\Box$
\vskip .3in

Recall that the Hilbert series  $\cg M\cd$ of  a finitely generated graded module $M$ over a homogeneous ring over a field is element of the ring ${\bf Z}[t,t^{-1},(1-t)^{-1}]\subset  
{\bf Z}\cg t\cd [t^{-1}]$. In general, any $S\in {\bf Z}[t,t^{-1},(1-t)^{-1}]$ can be written uniquely in the form
$$
S(t)=\sum_{i=0}^{D-1}(-1)^{i}e_{i}{{1}\over{(1-t)^{D-i}}}+F
$$
where $e_i\in {\bf Z}$ and $F\in {\bf Z}[t,t^{-1}]$. The coefficients $e_i$ can be computed as 
$e_{i}(M)={{\partial^{i}P}\over{i!}}(1)$, where $P(t)=S(t)(1-t)^{D}$. We call 
$$
Q(t)=\sum_{i=0}^{D-1}(-1)^{i}e_{i}{{t+D-1-i}\choose{D-1-i}}\in {\bf Q}[t]
$$
the polynomial associated to $S$. Its significance is that if we write $S=\sum_{i\in {\bf Z}}c_i t^i$, 
then $c_i=Q(i)$ for $i\gg 0$.

{\bf Remark 1.8.} In the case where $S(t)=\cg M\cd (t)$, we can take $D$ to be any integer $\geq \dim M$, 
and we define $e_i^D (M):=e_i$.
If $D=\dim M$, we simply set $e_i(M):=e_i^D(M)$. Notice that $e_0 (M)$ is the multiplicity (or degree)
of $M$. The polynomial associated to $\cg M\cd (t)$ is the Hilbert polynomial of $M$, which we denote by $[M](t)$.

\vskip .3in

{\bf Theorem 1.9.} {\sl  Write $e_{\l }(d_{1},\ldots ,d_{s})=
\sum_{{i_{1}\geq 1,\ldots ,i_{s}\geq
1}\atop{\scriptscriptstyle{i_{1}+\cdots +i_{s}=\l 
+s}}}\prod_{j=1}^{s}{{d_{j}}\choose{i_{j}}}$.

{\rm (a)} With the assumptions of Theorem 1.4(a),
$$
\eqalign{
e_{i}^{n-s}(I/\J )&=\sum_{k=0}^{i}e_{i-k}(d_{1},\ldots
,d_{s})e_{k}(R)-(-1)^{s-g}e_{s-g+i}(R/I)\cr
 &-(-1)^{s-g}\sum_{j=1}^{s-g}\sum_{k=0}^{s-g+i}(-1)^{j+k}\sum_{1\leq
i_{1}<\cdots <i_{g+j}\leq s}{{d_{i_{1}}+\cdots
+d_{i_{g+j}}+n-g-k}\choose{i+s-g-k}}\; e_{k}(\o /I^{j}\o )\cr}
$$
for $0\leq i\leq r-s$.

{\rm (b)} With the assumptions of Theorem 1.4(b),
$$
\eqalign{
e_{i}^{n-s}&(R/\K )=\sum_{k=0}^{i}e_{i-k}(d_{1},\ldots
,d_{s})e_{k}(R)\cr
 &+(-1)^{s-g}\sum_{j=1}^{s-g+1}\sum_{k=0}^{s-g+i}(-1)^{j+k}\sum_{1\leq
i_{1}<\cdots <i_{g+j-1}\leq s}{{d_{i_{1}}+\cdots
+d_{i_{g+j-1}}+n-g-k}\choose{i+s-g-k}}\; e_{k}(\o /I^{j}\o )\cr}
$$
for $0\leq i\leq r-s$.\medskip
}
\medskip

{\bf Proof.} We only prove part (a). First write
$$
\eqalign{
(-1)^{n-g}\sum_{j=1}^{s-g}(-1)^{j}\sum_{1\leq
i_{1}<\cdots <i_{g+j}\leq s}[\o /I^{j}\o ](-t+d_{i_{1}}+\cdots +d_{i_{g+j}})
&=\sum_{\l =0}^{n-g-1}(-1)^{\l}h_{\l}{{t+n-g-1-\l}\choose{n-g-1-\l}}\cr
&=\sum_{\l =0}^{n-s-1}(-1)^{\l}h_{\l}^{*}{{t+n-s-1-\l}\choose{n-s-1-\l}}\cr}
$$
and notice that $h_{i}^{*}=(-1)^{s-g}h_{i+s-g}$. Lemma 1.7 gives
$$
h_{i+s-g}=-\sum_{j=1}^{s-g}(-1)^{j}\sum_{1\leq
i_{1}<\cdots <i_{g+j}\leq
s}\sum_{k=0}^{i+s-g}{{d_{i_{1}}+\cdots
+d_{i_{g+j}}+n-g-k}\choose{i+s-g-k}}(-1)^{k}e_{k}(\o /I^{j}\o ).
$$

Now our assertion follows from Theorem 1.4(a) together with Lemma 1.6.

The proof of part (b) is similar, using Theorem 1.4(b) in place of Theorem 1.4(a).\hfill$\Box$
\medskip

{\bf Remark 1.10.} If in Theorem 1.9, $R$ is a polynomial ring in $n$
variables, then the formula in that theorem takes the following form:

$$
\eqalign{e_{i}^{n-s}(I/\J )&=e_{i}(d_{1},\ldots ,d_{s})-(-1)^{s-g}e_{s-g+i}(R/I)\cr
&-(-1)^{s-g}\sum_{j=1}^{s-g}\sum_{k=0}^{s-g+i}(-1)^{j+k}\!\!\!\sum_{1\leq
i_{1}<\cdots <i_{g+j}\leq s}{{d_{i_{1}}+\cdots
+d_{i_{g+j}}-g-k}\choose{i+s-g-k}}\; e_{k}(R/I^{j})\cr}
$$
for $0\leq i\leq r-s$;

$$
\eqalign{e_{i}^{n-s}(R/\K )&=e_{i}(d_{1},\ldots ,d_{s})\cr
&+(-1)^{s-g}\!\sum_{j=1}^{s-g+1}\sum_{k=0}^{s-g+i}(-1)^{j+k}\!\!\!
\!\!\!\sum_{1\leq
i_{1}<\cdots <i_{g+j-1}\leq s}{{d_{i_{1}}+\cdots
+d_{i_{g+j-1}}-g-k}\choose{i+s-g-k}}\; e_{k}(R/I^{j})\cr}
$$
for $0\leq i\leq r-s$.\medskip

{\bf Proof.}  Notice that $\o /I^{j}\o\simeq R/I^{j} (-n)$ and proceed as in the proof of Theorem 1.9.\hfill$\Box$
\medskip

{\bf Corollary 1.11.} Let $R$ be a homogeneous ring over a field, write
$n=\dim R$, $\o =\o_{R}$, and assume that $R$ is Gorenstein locally in
\tc $r\geq s$. Let $I$ be a homogeneous ideal of height $g$
satisfying $G_{s}$, let $f_{1},\ldots ,f_{s}$ be forms contained in $I$
of degrees $d_{1},\ldots ,d_{s}$, write $\J =(f_{1},\ldots ,f_{s})$,
$\K =\J :I$, and assume that locally in \tc $r$, $\hbox{depth}\,
R/I^{j} \geq \dim R/I-j+1$ for $1\leq j\leq s-g$.

Then ${\rm ht}\, \K \geq r+1$ if and only if ${\rm ht}\, \K \geq s$ and
$$
\eqalign{
(-1)^{s-g}e_{0}(R)\prod_{j=1}^{s}d_{j}&= e_{s-g}(R/I)\cr
& +\sum_{j=1}^{s-g}
\sum_{k=0}^{s-g} (-1)^{j+k}\sum_{1\leq i_{1}<\cdots <i_{g+j}\leq
s}{{d_{i_{1}}+\cdots +d_{i_{g+j}}+n-g-k}\choose{s-g-k}}
e_{k}(\o /I^{j}\o ).\cr}
$$

{\bf Proof.} One uses Theorem 1.9(a) and [20, 1.7(a)].\hfill$\Box$ 
\medskip

\bigskip

\goodbreak{\bigrm 2. Hilbert series of powers of
ideals and degrees of residual intersections.}\medskip 

{\bf 2.1 Computing Hilbert series of powers.}\medskip

Motivated by Remark 1.5, showing the usefulness of the
Hilbert series of the powers of an ideal, we will focus
here on the following question: to what extent does the Hilbert
series of the first powers of an ideal determine the Hilbert series of the next
powers? 

The following Lemma tells us a useful property of a general set of elements of an ideal.

{\bf Lemma 2.1.} {\sl Let $R$ be a standard graded Cohen-Macaulay ring
over an infinite field $k$, let $I$ a homogeneous ideal, generated by forms of degrees at most $d$,
and let $r$ be an integer with $0\leq r\leq\dim R$ . Further assume that
$I$ satisfies   $G_{r+1}$ and has sliding depth locally in codimension $r$.

Given $d_{i}\geq d$ for $1\leq i\leq r+1]$, there exists a Zariski dense open
subset $\Omega$ of the affine 
$k$-space $I_{d_{1}}\times\cdots\times I_{d_{r+1}}\simeq {\bf
A}_{k}^{N}$ (where $N=\sum_{i=1}^{r+1}\dim_{k}I_{d_{i}}$) such that if
$(f_{1},\ldots ,f_{r+1})\in \Omega$~:

{\rm (a)} The ideal $(f_{1},\ldots ,f_{r+1})$ coincides with $I$ locally
up to codimension $r$.

{\rm (b)} $f_{1},\ldots ,f_{r+1}$ is a $d$-sequence locally up to
codimension 
$r$.}
\medskip

{\bf Proof.} For part (a) we refer to [2, 1.6 (a)] and [10, 3.9], while
(b) follows from [6, 3.6(b)].\hfill $\Box$
\medskip

\def\ra{\rightarrow}

We now choose some degrees (e.g. $d_{i}=d$ for all $i$) and
polynomials $f_{i}$ as above. Part two of the lemma implies that the
approximation complexes corresponding to $f_{1},\ldots ,f_{r+1}$, and
therefore their different graded components,  
$$
{\cal M}_{p}:\quad 0\ra H_{r+1}\otimes S_{p-r-1}\ra H_{r}\otimes
S_{p-r}\ra \cdots \ra H_{0}\otimes S_{p}\ra 0
$$
have  homology in positive degrees supported in codimension at least
$r+1$. Moreover, $H_{0}({\cal M}_{p})$ coincides with 
$I^{p}/I^{p+1}$ locally up to codimension $r$. See [12] for all of
these facts.\medskip 

Recall that in this sequence, $H_{q}$ stands for the $q$-th homology
module of the Koszul complex ${\bf K }(f_{1},\ldots ,f_{r+1};R)$ and
$S_{q}$ is the free $R$-module generated by monomials of degree $q$ in
$r+1$ variables. The maps are homogeneous of degree 0 in the graded case,
with the usual weights on the Koszul complex and similarely the weight of
$s_{i_{1}}\cdots s_{i_{q}}\in S_{q}$ is $\deg (f_{i_{1}})+\cdots +\deg
(f_{i_{q}})$. 

\def\cg{\lbrack\!\lbrack}
\def\cd{\rbrack\!\rbrack}
\def\ra{\rightarrow}
\def\lra{\longrightarrow}
\def\o{\omega}
\def\s{\sigma}

\def\d{\delta}

Therfore we have the following equalities~:
$$
\cg I^{p}/I^{p+1}\cd (t)\egr \sum_{i=0}^{p}(-1)^{i}s_{p-i}(t^{d_{1}},\ldots
,t^{d_{r+1}})\cg H_{i}\cd (t)\leqno{(1)_{p}}
$$
for every $p\geq 0$, where $s_{j}$ stands for the sum of all the
monomials of degree $j$ in $r+1$ variables (complete symmetric
functions). These express the Hilbert series of  
the modules $H_{0},\ldots ,H_{p}$ in terms of the ones of $R/I,\ldots ,
I^{p}/I^{p+1}$, and vice versa.

If $I$ has height $g$, then $H_{q}=0$ for $q>r+1-g$. We
therefore immediately see that the Hilbert series of all the modules
$I^{p}/I^{p+1}$ are determined, up to $r$-equivalence, by the
knowledge of the Hilbert series of $I^{p}/I^{p+1}$ for $0\leq p\leq
r+1-g$, up to $r$-equivalence.\medskip

We now assume that, in addition, $R$ is Gorenstein with $a$-invariant $a:=a(R)$
and is strongly Cohen-Macaulay locally in codimension  $r$. In this case we can use
the self-duality of the homology of the Koszul complex to see that 
only half of the information about the Koszul homology is needed. Indeed, the
structure of graded alternating algebra on the homology of the Koszul
complex gives a graded map of degree 0,
$$
H_{p}\lra {\rm Hom}_{R/I}(H_{r+1-g-p},H_{r+1-g}),
$$
which is an isomorphism up to codimension $r$ by a theorem of Herzog
(see [11, 2.4.1]).

Now we have a collection of graded maps of degree 0 that connect the
following modules,
$$
\eqalign{
{\rm Hom}_{R/I}(H_{r+1-g-p},H_{r+1-g})
&\eqr {\rm Hom}_{R/I}(H_{r+1-g-p},{\rm
Ext}^{g}_{R}(R/I,R)[-(d_{1}+\cdots +d_{r+1})])\cr
&\eqr {\rm
Hom}_{R/I}(H_{r+1-g-p},\o_{R/I}[-a-(d_{1}+\cdots +d_{r+1})])\cr
&\eqr \o_{H_{r+1-g-p}}[a+(d_{1}+\cdots +d_{r+1})].\cr}
$$
From the Cohen-Macaulayness of the modules $H_{p}$, locally in
codimension at most $r$, we therefore have:
$$
\cg H_{p}\cd (t)\egr t^{a+(d_{1}+\cdots +d_{r+1})}(-1)^{\dim R/I}\cg
H_{r+1-g-p}\cd (t^{-1}). \leqno{(2)_{p}}
$$

Moreover, the Euler characteristic of the homology of the
Koszul complex depends only upon the degrees 
$d_1,\ldots ,d_{r+1}$, namely
$$
\sum_{p=0}^{r+1-g}(-1)^{p}\cg H_{p}\cd =
\Delta_{r+1}\cg R\cd \egr 0.\leqno{(3)}
$$ 

We now have put together all the formulas needed to effectively
compute what we state in the next theorem.\medskip

Remark that the result of the computation does not depend on the choice
of the $d_{i}$'s. Thus we may choose $d_{i}=0$ for all $i$, the
intermediate steps have no meaning (e.g. $\cg H_{p}\cd$ may not have
positive coefficients), but the information that we extract from the
computation is the same. Out of this remark, one may use the following~:

$$
\cg I^{p}/I^{p+1}\cd (t)\egr
\sum_{i=0}^{p}(-1)^{i}{{r+p-i}\choose{r}}\cg H_{i}\cd
(t), \leqno{(1)_{p}} 
$$
$$
\cg H_{p}\cd (t)\egr t^{a}(-1)^{\dim R/I}\cg
H_{r+1-g-p}\cd (t^{-1}), \leqno{(2)_{p}}
$$
$$
\sum_{p=0}^{r+1-g}(-1)^{p}\cg H_{p}\cd \egr 0.\leqno{(3)}
$$ 

{\bf Theorem 2.2. } {\sl Let $R$ be a homogeneous algebra, $I$
an homogeneous $R$-ideal, and let us suppose that, locally in
\tc  $r$, $I$ is licci purely of \tc $g$ and satisfies
$G_{r+1}$. Given the
Hilbert series of $I^{p}/I^{p+1}$ for $0\leq p\leq
\left[ {r-g}\over{2} \right]$, up to $r$-equivalence, the
Hilbert series of $I^{p}/I^{p+1}$ can be computed for all $p$,
up to $r$-equivalence, by the formulas above.}\medskip
 
{\bf Proof.} Let us choose some sufficiently big integers $d_{1},\ldots
,d_{r+1}$ (one can treat them as unknowns, or take all of them equal to some
fixed number or unknown $d$). Let us put $q=r+1-g$.

First, using $(1)_{p}$ for $0\leq p\leq \left[ {q-1}\over{2} \right] $,
we get the Hilbert series of $H_{p}$, for $p$ in the same
range, up to \tc $>r$ terms. 

Then, from $(2)_{p}$ for the same $p$'s, we get the series of
$H_{q},\ldots ,H_{q-\left[{q-1}\over{2} \right]}$. Therefore, we get the
series of all the $H_{p}$'s (up to $r$-equivalence) if $q$ is
odd~; and all of them but one, namely $H_{q/2}$, if $q$ is even. 
In case $q$ is even, we get the series of $H_{q/2}$ using
$(3)$. 

Therefore we know the series of all the modules $H_{p}$ in every
case, and can use $(1)_{p}$ to get the ones of $I^{p}/I^{p+1}$ for any
$p$. 

Notice that the $d_{i}$'s appear in every formula, but should disappear
at the end~!
\hfill $\Box$
\medskip

{\bf Example 2.3.} {\sl If $X\subseteq {\Bbb P}^{n}$ is an equidimensionnal
locally complete intersection scheme, and ${\cal I}_{X}$ the
corresponding ideal sheaf, the Hilbert polynomials of the sheaves ${\cal
I}_{X}^{p}/{\cal I}_{X}^{p+1}$ are all determined by the ones for $0\leq
p\leq \left[ {\dim X}\over{2} \right]$.}\medskip

{\bf Proof.} It is the case where $r=n$ and $g=n-\dim X$.
\hfill $\Box$
\medskip

{\bf Example 2.4.} If $X\subseteq {\Bbb P} ^{n}$ is an equidimensionnal
locally complete intersection threefold, and ${\cal I}_{X}$ the
corresponding ideal sheaf, the Hilbert polynomial of the sheaves ${\cal
I}_{X}^{p}/{\cal I}_{X}^{p+1}$ are all determined by the Hilbert polynomials of
$X$ and the one of the conormal bundle ${\cal I}_{X}/{\cal
I}_{X}^{2}$. Moreover the  coefficients of terms in degrees 3 and 2 in
the Hilbert polynomial of the conormal bundle are determined by those
of the Hilbert polynomial for $X$.\medskip

Informally, we have the following picture for the determination of the
highest $r-g+1$ coefficients of the Hilbert polynomial of the powers of
an ideal $I$ of codimension $g$ that is locally complete intersection up to
codimension $r$~:

$$
\matrix{
\hbox{\rm HP coef.}& &\cdot&\cdot&g&\cdot&\cdot&\cdot&\cdot&r&\cdot&\cdot& \cr
                 &   &  &   &     &    &    &    &    &   &  & & \cr
         R/I &\cdots&  0& 0    &\blacksquare   &\blacksquare  &\blacksquare  &\blacksquare  &\blacksquare  &\blacksquare & ?&?&\cdots\cr
      I/I^{2} &\cdots& 0 & 0    &\Box &\Box&\blacksquare  &\blacksquare  &\blacksquare  & \blacksquare& ?&?&\cdots \cr
    I^{2}/I^{3}&\cdots& 0& 0    &\Box &\Box&\Box&\Box&\blacksquare  & \blacksquare& ?&?&\cdots \cr
  I^{3}/I^{4} &\cdots& 0& 0    &\Box &\Box&\Box&\Box&\Box&\Box& ?&?&\cdots \cr
 I^{4}/I^{5} &\cdots& 0& 0    &\Box &\Box&\Box&\Box&\Box&\Box&  ?&?&\cdots\cr
\vdots&&\vdots&\vdots&\vdots&\vdots&\vdots&\vdots&\vdots&\vdots&\vdots&\vdots&
\cr}
$$

$\blacksquare$~: needed as input.\par
$\Box$~: may be computed from the others.\par
$?$~: not concerned.\medskip

\bigskip

\goodbreak{\bf 2.2 Hilbert polynomials of powers of an ideal}\medskip

We will treat the example of an equidimensional locally complete
intersection threefold in projective $n$-space. 

Let us abreviate $e_{i}(p)=e_{i}(I^{p}/I^{p+1})$. Theorem 2.2 asserts
that all the coefficients $e_{i}(p)$ are determined by six of them.

With the help of a computer algebra system, one gets the following
formulas.

{\bf Formulas 2.5}
$$
\eqalign{
e_{0}(p)&={{g+p-1}\choose{p}} e_{0}(0),\cr
e_{1}(p)&=g{{g+p-1}\choose{p-1}} e_{0}(0)
+{{(g+2p)}\over{(g+p)}}{{g+p}\choose{p}}e_{1}(0),\cr
e_{2}(p)&={{g(g+1)}\over{2}}{{g+p-1}\choose{p-2}}e_{0}(0)
+(g+1){{g+p-2}\choose{p-2}}e_{1}(0)\cr
&-{{(p-1)g}\over{(g+p)}}{{g+p+1}\choose{p}}
e_{2}(0)+{{g+p}\choose{p-1}}e_{2}(1),\cr
e_{3}(p)&={{g(g+1)(g+2)}\over{6}}{{g+p-1}\choose{p-3}}e_{0}(0)
-{{(g+1)(g+2)}\over{2}}{{g+p-1}\choose{p-2}}e_{1}(0)\cr
&-{{p(p-1)(g+2)}\over{(g+p)}}{{g+p+1}\choose{p-2}}e_{2}(0)
+(g+2){{g+p}\choose{p-2}}e_{2}(1)\cr
&+{{(p-1)(g+2p)}\over{(g+p)}}{{g+p+1}\choose{p}}e_{3}(0)
+{{(g+2p)}\over{(g+2)}}{{g+p}\choose{p-1}}e_{3}(1).\cr}
$$

Notice that these formulas remains valid for the case of any scheme
$X\subseteq {\Bbb P}^{n}$ that is locally a complete intersection in
dimension $\geq \dim X-3$.\medskip

As a first guess, one may hope that, at least with
some strong hypotheses on $X$, the Hilbert polynomial of the powers are
determined by the Hilbert polynomial of $X$. This is even not true for
complete intersections, due to the following computation. 

Suppose that $X$ is a global complete intersection of codimension $g$
and denote by $\s_{1},\ldots ,\s_{g}$ the symmetric functions on the
degrees of the defining equations of $X$. Setting
$$
\a_{1}=\s_{1}-g,\quad \a_{2}=\s_{1}^{2}-2\s_{2}-g,\quad 
\a_{3}=\s_{1}^{3}-3\s_{1}\s_{2}+3\s_{3}-g,
$$ 
one gets,
$$
\eqalign{
e_{0}(0)&=\s_{g},\cr
e_{1}(0)&={{\s_{g}}\over{2}}\a_{1},\cr
e_{2}(0)&={{\s_{g}}\over{24}}(3\a_{1}^{2}-6\a_{1}+\a_{2}),\cr
e_{3}(0)&={{\s_{g}}\over{48}}(\a_{1}^{3}-6\a_{1}^{2}+8\a_{1}+\a_{1}\a_{2}-2\a_{2}).\cr}
$$
Notice that these formulas imply that $e_{3}(0)$, the fourth coefficient of the
Hilbert polynomial, is a rational function of the first three :

{\bf Remark 2.6} If $e_i$ denotes the $i$-th coefficient of the Hilbert polynomial of a global complete intersection 
of dimension at least 3 in
a projective space, then
$$
e_3=e_2-{{e_1 e_2}\over{e_0}}+{{e_1}\over{6}}-{{e_1^2}\over{2e_0}}+{{e_1^3}\over{3e_0^2}}.
$$

Now, using the expansion
$$
t^{d_{1}}+\cdots
+t^{d_{g}}=g+(g+\a_{1})(t-1)+(\a_{2}-\a_{1}){{(t-1)^{2}}\over{2}}
+(\a_{3}-3\a_{2}+2\a_{1}){{(t-1)^{3}}\over{6}}+\cdots ,
$$
one can compute the coefficients $e_{i}(1)$ for $1\leq i\leq 3$. 
The only place where $\a_{3}$ appears is in
$e_{3}(1)={{\s_{g}}\over{6}}\a_{3}+\cdots$.\medskip

If one chooses two collections of degrees such that the first, second
and $4$-th symmetric functions are equal but the third one differs, one
gets an example of two complete intersections of dimesion three in $\PP^{7}$ having the same Hilbert
polynomials (but distinct Hilbert functions~!) such that the constant
term of the Hilbert polynomials of their conormal bundle are distinct.
Such examples were given to us by Benjamin de Weger, the two ``smallest''
ones are (1,6,7,22)-(2,2,11,21) and (2,6,7,15)-(3,3,10,14). He also gave
an infinite collection of them, and Noam Elkies gave a rational 
parametrization of all the solutions (after a linear change of
coordinates the solutions are parametrized by a quadric in ${\bf
P}^{5}$).

\goodbreak

{\bf 2.3 The degree of the residual}\bigskip

Let us suppose that the projective scheme $X$ of dimension $D$ is
locally a complete 
intersection in codimension at most $s$, and use our formulas and the
above computations to derive the degree of the codimension $s$ part of
an $s$-residual intersection.\medskip

We will treat the cases where $\d =s-g$ is less or equal to three. As
before $\s_{i}$ stands for the $i$-th symmetric function on
$d_{1},\ldots ,d_{s}$. We will also set, to simplify some formulas,
$e'_{i}(p)=\sum_{j=0}^{p}e_{i}(j)$, which is the $i$-th coefficient of the
Hilbert polynomial $R/I_{X}^{p}$.\medskip

$\bullet$ If $\d =0$, $e_{0}^{D}(R/\K )=\s_{s}-e_{0}(0)$
(B{\'e}zout).\medskip

$\bullet$ If $\d =1$, $e_{0}^{D-1}(R/\K )=\s_{s}-(\s_{1}-g)e_{0}(0)+2e_{1}(0)$, as
proved by St{\"u}ckrad in [18].\medskip

$\bullet$ If $\d =2$, using the Taylor  expansion of $\s_{1}(t^{d_{1}},\ldots
,t^{d_{s}})$,
$$
\eqalign{
\s_{1}(t^{d_{1}},\ldots
,t^{d_{s}})&=s+\s_{1}(t-1)+(\s_{1}^{2}-\s_{1}-2\s_{2}){{(t-1)^{2}}\over{2}}\cr
&+(\s_{1}^{3}-3\s_{1}^{2}+2\s_{1}-3\s_{2}(\s_{1}-2)+3\s_{3}){{(t-1)^{3}}\over{6}}+\cdots
,\cr}
$$
one recovers the formula given by Huneke and
Martin in [16],
$$
e_{0}^{D-2}(R/\K )=\s_{s}-\left( \s_{3}-g\s_{2}+{{g+1}\choose{1}}\right)
e_{0}(0)+(2\s_{1}-(g+1))e_{1}(0)+(g+1)e_{2}(0)-e'_{2}(1).
$$

$\bullet$ If $\d =3$, using the following Taylor expansion,
$$
\eqalign{
\s_{2}(t^{d_{1}},\ldots ,t^{d_{s}})&={{s(s-1)}\over{2}}+(s-1)\s_{1}(t-1)
+((s-1)(\s_{1}^{2}-\s_{1}-2\s_{2})+2\s_{2}){{(t-1)^{2}}\over{2}}\cr
&+((s-1)(\s_{1}^{3}-3\s_{1}^{2}+2\s_{1})-3(s-2)\s_{2}(\s_{1}-2)+3(s-4)\s_{3})
{{(t-1)^{3}}\over{6}}+\cdots}
$$
one gets from Theorem 1.4,
$$
\eqalign{
e_{0}^{D-3}(R/\K )&=\s_{s}-\left(
\s_{3}-g\s_{2}+{{g+1}\choose{2}}\s_{1}-{{g+2}\choose{3}}\right)
e_{0}(0)-(2\s_{2}-(g+1)\s_{1})e_{1}(0)\cr
&+((g+1)\s_{1}-(g+2)(g+3))e_{2}(0)-(\s_{1}-(g+2))e'_{2}(1)+2(g+1)e_{3}(0)
+2e'_{3}(1).\cr}
$$
 
\bigskip

\goodbreak{\bigrm 3. Applications to secant varieties}\bigskip

{\bf Theorem 3.1.} {\it Let $k$ be a perfect field, $X\subset {\bf P}^{N}_{k}$ an 
eqidimensional
subscheme of dimension two with at most isolated licci Gorenstein singularities, $A$ 
its homogeneous coordinate ring, $\om :=\om_{A}$ the canonical module, and 
$\Om :=\Om_{A/k}$ the module of differentials. \smallskip
{\rm (a)} One has
$$
e_{0}(A)^{2}+14e_{0}(A)-16e_{1}(A)+4e_{2}(A)\geq e_{2}(\om \otimes_{A}\om )+e_{2}(\Om ).
$$

{\rm (b)} In case the singularities of $X$ have embedding codimension at most two, then
equality holds in {\rm (a)} if and only if the secant variety of $X$ is deficient, i.e.
$$
\dim \hbox{\rm Sec}(X)<5.
$$
}
{\it Proof.} We may assume that $k$ is infinite. We define the ring $R$ and the $R$-ideal 
$I$ via the exact sequence
$$
0\lra I\lra R:=A\otimes_{k}A\buildrel{\rm mult}\over{\lra} A\lra 0.
$$
Recall that $\Omega \simeq I/I^2$. The ring
$R$ is a standard graded $k$-algebra of 
dimension 6 with $\om_{R}=\om\otimes_{k}\om$ and $\codim NG(R)>2$. The ideal $I$ has height 3, and
is generated by linear forms. Moreover, $I$ satisfies $G_{5}$ and, in the setting of (b),
even $G_{6}$. Indeed,
for any $\ip \in V(I)$, one has $\mu (I_{\ip})=\mu 
(\Om_{\ip})\leq \hbox{\rm ecodim}(A_{\ip})+\dim A\leq \dim R_{\ip}$
if $\dim R_{\ip}\leq 4$ or, in the setting of (b), $\dim R_{\ip}\leq 5.$
In addition, for every $\ip \in V(I)$ with $\dim R_{\ip}\leq 5$, we have 
$\depth (I/I^{2})_{\ip}=\depth \Om_{\ip} \geq \dim A_{\ip}-1$.
To see this, we write $A_{\ip}\simeq S/\JJ$ with $S$ a regular
local ring and $\JJ$ a licci Gorenstein ideal. This is possible because $A_\ip$ is
licci and Gorenstein. By [4, 6.2.11 and 6.2.12]
the module $\JJ /\JJ^2$ is Cohen-Macaulay. Thus the natural complex
$$
0\to \JJ /\JJ^2 \to \Omega_{S/k} \otimes_S A_\ip \simeq \oplus A_\ip 
\to 
\Omega_{A_\ip/k}
\simeq
\Omega_\ip
\to
0
$$
is exact and shows that depth $\Omega_\ip \geq {\rm dim}\ A_\ip -1$.

Now let $\ia$ be an $R$-ideal generated by 5 general linear forms in $I$. 
Notice that the five general linear forms in $I$ that generate $\J$ are a 
weak $R$-regular sequence off $V(I)$, hence off $V(\J)$.
Since $\codim NG(R)>0$ and 
$\codim NG(R)\cap V(I)>5$ it follows that $\codim NG(R)\cap V(\ia )>5$. 
By [2, 1.4] one has 
$\he (\ia :I)\geq 5$ as $I$ satisfies $G_{5}$, and $\he (I+(\ia :I))\geq 6$ in (b) 
as $I$ is $G_{6}$. Thus in the setting of (b), [19] shows that $\he (\ia :I)\geq 6$ 
if and only if  the analytic spread $\ell (I)$ is at most 5. On the other hand
$\dim \hbox{\rm Sec}(X)=\ell (I)-1$ according to [17]. Hence $\dim \hbox{\rm Sec}(X)<5$
 if and only if $e^{1}_{0}(I/\ia )=0$.

We now apply Theorem 1.9(a) with $r=s=5$ and $g=3$. The theorem yields
$$
\leqno{(1)}\quad  e^{1}_{0}(I/\ia )=e_{0}(R)-e_2(A)-\sum_{j=1}^{2}\sum_{k=0}^{2}(-1)^{j+k}{{5}\choose{j+3}}
{{6+j-k}\choose{2-k}}e_{k}(\om_{R}/I^{j}\om_{R}).
$$
Thus the present theorem follows once we have shown that the right hand side of (1)
equals 
$$
\leqno{(2)}\quad  e_{0}(A)^{2}+14e_{0}(A)-16e_{1}(A)+4e_{2}(A)-e_{2}(\om \otimes_{A}\om )-e_{2}(\Om ).
$$

From [5, IX 2.1] one obtains the isomorphisms of $R$-modules 
$$
\om_{R}/I\om_{R}\cong \om_{R}\otimes_{R}R/I\cong (\om \otimes_{k}\om )
\otimes_{A\otimes_{k}A}A
\cong \om\otimes_{A}(\om\otimes_{A}A ) \cong \om\otimes_{A}\om
$$
and, since $\codim NG(R)\cap V(I)>5$,
$$
I\om_{R}/I^{2}\om_{R}\isof \om_{R}\otimes_{R}I/I^{2}\cong 
(\om \otimes_{k} \om )\otimes_{A\otimes_{k}A}\Om\cong
\om\otimes_{A}(\om \otimes_{A}\Om )\cong 
(\om \otimes_{A}\om)\otimes_{A}\Om .
$$
Therefore the right hand side of (1) becomes
$$
\leqno{(3)}\quad  e_{0}(A)^{2}-7e_{0}(A)-e_{2}(A)-23e_{1}(\om^{\otimes 2})
+4e_{2}(\om^{\otimes 2})+7e_{1}(\om^{\otimes 2}\otimes \Om)
-e_{2}(\om^{\otimes 2}\otimes \Om).
$$
Here and in what follows tensor products are taken over the ring $A$.

We are now going to express the Hilbert coefficients $e_{1}(\om^{\otimes 2})$, 
$e_{1}(\om^{\otimes 2}\otimes \Om)$ and $e_{2}(\om^{\otimes 2}\otimes \Om)$ in 
terms of the Hilbert coefficients of $A$, $e_{2}(\om^{\otimes 2})$ and $e_{2}(\Om )$.
First notice that for any finitely generated graded $A$-module $M$,
$$
\leqno{(4)}\quad  e_{i}(M(-1))=e_{i}(M)+e_{i-1}(M).
$$
Hence, by Lemma 1.7 and Remark 1.8, 
$$
\leqno{(5)}\quad  e_{1}(\om )=3e_{0}(A)-e_{1}(A)\quad {\rm and}
\quad e_{2}(\om )=3e_{0}(A)-2e_{1}(A)+e_{2}(A).
$$

Since $\om$ is free of rank 1 locally in codimension 1, there is a complex of graded $A$-modules
$$
\leqno{(6)}\quad  0\lra Z\lra A(-a)^{2}\lra \om\lra 0
$$
for some $a\gg 0$, that is exact in codimension 1. It induces complexes
$$
\leqno{(7)}\quad  0\lra Z\otimes A(-ja+a)^{j}\lra A(-ja)^{j+1}\lra \sym_{j}(\om )
\isoo \om^{\otimes j}
\lra 0
$$
that are likewise exact locally in codimension 1. Now (6) yields $e_{1}(Z)
=2e_{1}(A(-a))-e_{1}(\om )$ and 
then (4), (5) and (7) show that for every $j\geq 0$,
$$
\leqno{(8)}\quad  e_{1}(\om^{\otimes j})=3je_{0}(A)-(2j-1)e_{1}(A).
$$

Now, we treat the first Hilbert coefficient of $\om^{\otimes 2}\otimes\Om$. Since $\om$ is free of rank 1 locally 
in codimension 2, we also have $\om^{*}\build\lra_{2}^{\sim}
\hom_{A}(\om^{\otimes 2},\om )$, which by Lemma 1.7 and Remark 1.8 gives
$$
\leqno{(9)}\quad  e_{1}(\om^{*})=3e_{0}(\om^{\otimes 2} )-e_{1}(\om^{\otimes 2} )\ \hbox{and}\ 
e_{2}(\om^{*})=3e_{0}(\om^{\otimes 2} )-2e_{1}(\om^{\otimes 2} )+e_{2}(\om^{\otimes 2} ).
$$
Furthermore,
$$
\leqno{(10)}\quad \om^{\otimes 2}\otimes \om^{*}\build\lra_{2}^{\sim}
\hom_{A}(\om ,\om^{\otimes 2})\build\longleftarrow_{2}^{\sim} \om .
$$

As $\Om$ is free of rank 3 locally in codimension 1, there is an exact sequence of
graded $A$-modules
$$
0\lra A(-1)^{2}\lra \Om \lra C\lra 0
$$
where $C$ is free of rank 1 locally in codimension 1. Thus 
$$
C\otimes \bigwedge^{2}(A(-1)^{2})\build\lra_{1}^{\sim}\bigwedge^{3}\Om 
\build\lra_{1}^{\sim}\om ,
$$
which gives a complex
$$
0\lra A(-1)^{2}\lra \Om \lra \om (2)\lra 0
$$
that is exact in codimension 1. Tensoring with $\om^{\otimes 2}$ we obtain
$$
0\lra \om^{\otimes 2}(-1)^{2}\lra \om^{\otimes 2}\otimes\Om \lra \om^{\otimes 3} (2)\lra 0.
$$
Since this complex is exact in codimension 1, (4) and (8) imply that
$$
\leqno{(13)}\quad  e_{1}(\om^{\otimes 2}\otimes\Om )=21e_{0}(A)-11e_{1}(A).
$$

Next, we turn to the second Hilbert coefficient of $\om^{\otimes 2}\otimes\Om$. 
Write $\hbox{---}^{*}:=\hom_{A}(\hbox{---},A)$ and $e:=N-2$. Increasing $N$ if needed, we may assume $e\geq 2$.
We define a graded $A$-module $E$ via the exact sequence
$$
\leqno{(14)}\quad  0\lra E\lra A(-1)^{e+3} \lra \Om\lra 0.
$$
Notice that $E$ has rank $e$, and is free locally in 
codimension 1 and Cohen-Macaulay locally in codimension 2. Futhermore (14) gives
$$
\leqno{(15)}\quad ( \bigwedge^{e}E) ^{**}
\build\lra_{2}^{\sim}( \bigwedge^{3}\Om
) ^{*}(-e-3)\build\longleftarrow_{2}^{\sim}\om^{*}(-e-3).
$$

As $E^{*}$ is free locally in codimension $1$ and $\rk E^{*}-1\geq 1$, 
there exists a 
homogeneous element $f\in E^{*}$ of degree $c\gg 0$ whose order ideal $(E^{*})^{*}(f)$
has height at least 2 (see [8]). However, the ideals $E^{**}(f)$ and $J:=f(E)$
coincide locally in codimension 1 since $E$ is reflexive locally in codimension 1. Hence
$\he J=\he E^{**}(f)\geq 2$.
The map $f$ induces an exact sequence of graded $A$-modules
$$
0\lra E_{e-1}\lra E_{e}:=E\lra J_{e}(c_{e}):=J(c)\lra 0.
$$
Repeating this procedure, if needed, we obtain a filtration
$$
\leqno{(16)}\quad E_{1}\subset E_{2}\subset\cdots\subset E_{e}\quad\hbox{with}\quad 
E_{i}/E_{i-1}\cong J_{i}(c_{i}),
$$
where $J_{i}$ are homogeneous $A$-ideals of height at least 2.
Thus  $E_{i}$ has rank $i$, is free in codimension 1 and Cohen-Macaulay 
in codimension 2, and 
$$
( \bigwedge^{i-1}E_{i-1})^{**} (c_{i})
\build\longleftarrow_{2}^{\sim}
( \bigwedge^{i-1}E_{i-1}) ^{**}\otimes(J_{i}(c_{i}))^{**}
\build\lra_{2}^{\sim}
( \bigwedge^{i}E_{i})^{**}.
$$
Since $E_1$ is reflexive locally in codimension $2$, it follows that
$$
E_{1} \build\lra_{2}^{\sim} E_{1}^{**}\build\cong_{2}^{}
 ( \bigwedge^{e}E)^{**}(-\sum_{i=2}^{e}c_{i}),
$$
which together with (15) implies
$$
\leqno{(17)}\quad E_{1}\build\cong_{2}^{}\om^{*}(-e-3-\sum_{i=2}^{e}c_{i}).
$$

The exact sequence
$$
0\lra J_{i}(c_{i})\lra A(c_{i})\lra (A/J_{i})(c_{i})\lra 0
$$
yields a complex 
$$
0\lra \om^{\otimes 2}\otimes J_{i}(c_{i})\lra \om^{\otimes 2}(c_{i})\lra (\om^{\otimes 2}
/\om^{\otimes 2}J_{i})(c_{i})\lra 0
$$
that is exact in codimension 2. Since $\he J_{i}\geq 2$ and $\om^{\otimes 2}$ 
is free of rank 1 locally in codimension 2, it follows that $e^{1}_{0}((A/J_{i})
(c_{i}))=e^{1}_{0}((\om^{\otimes 2}/\om^{\otimes 2}J_{i})(c_{i}))$.
We conclude
$$
\leqno{(18)}\quad e_{2}(J_{i}(c_{i}))-e_{2}(\om^{\otimes 2}\otimes J_{i}(c_{i}))
=e_{2}(A(c_{i}))-e_{2}(\om^{\otimes 2}(c_{i})).
$$

Tensoring (14) and (16) with $\om^{\otimes 2}$ and using (18) and (17) we obtain
$$
\leqno{(19)}\quad 
\eqalign{
e_{2}(\om^{\otimes 2}\otimes \Om )-e_{2}(\Om )&=e_{2}(\om^{\otimes 2}(-1)^{e+3})-
e_{2}(A(-1)^{e+3})+\sum_{i=2}^{e}(e_{2}(A(c_{i}))-e_{2}(\om^{\otimes 2}(c_{i})))\cr
&+e_{2}(\om^{*}(-e-3-\sum_{i=2}^{e}c_{i}))-e_{2}(\om^{2}\otimes 
\om^{*}(-e-3-\sum_{i=2}^{e}c_{i})).\cr}
$$

Combining (19) with (4), (8),  (9), (10), (5), we deduce
$$
\leqno{(20)}\quad  e_{2}(\om^{\otimes 2}\otimes \Om )=-12e_{0}(A)+8e_{1}(A)-5e_{2}(A)
+5e_{2}(\om^{\otimes 2})+e_{2}(\Om ).
$$
Substituting (8), (13), (20) into (3), we conclude that (3) and (2) coincide.\hfill$\Box$
\bigskip

{\bf Corollary 3.2.} Let $k$ be a perfect field, $X\subset {\bf P}^{4}_{k}$ an equidimensional
subscheme of dimension two with at most isolated Gorenstein singularities and $A$ its 
homogeneous coordinate ring. One has
$$
e_{0}(A)^{2}+14e_{0}(A)-16e_{1}(A)+4e_{2}(A)=e_{2}(\om_{A} \otimes_{A}\om_{A} )
+e_{2}(\Om_{A/k} ).
$$
\medskip

{\bf Remark.} The inequality in Theorem 3.1 can be replaced by
$$
e_{0}(A)^{2}+5e_{0}(A)-10e_{1}(A)+4e_{2}(A)\geq e_{2}(\om_{A}^{*})
+e_{2}(\Om_{A/k} ).
$$

{\it Proof.} Use equalities (9) and (8)  in the proof of Theorem 3.1.\hfill$\Box$
\bigskip

{\bf Corollary 3.3.} Let $k$ be a field, $X\subset {\bf P}^{N}_{k}$ an equidimensional
smooth
subscheme of dimension two, $H$ the class of the hyperplane section, $K$ the canonical 
divisor, and $c_{2}$ the second Chern class of the cotangent bundle of $X$. One has
$$
(H^{2})^{2}\geq 10H^{2}+5HK+K^{2}-c_{2},
$$
and equality holds if and only if $\dim \hbox{\rm Sec}(X)<5$.
\medskip
{\it Proof.} The Riemannn-Roch theorem in dimension two gives
$$
\chi (X,E)={{1}\over{2}}[c_{1}(E)^{2}-2c_{2}(E)-c_{1}(E)K_{X}]
+\rk E\chi (X,{\cal O}_{X}). 
$$
If $D$ is a divisor  this equality specializes to
$$
\eqalign{
\chi (D+nH)&={{1}\over{2}}H^{2}n^{2}+(DH-{{1}\over{2}}KH)n+{{1}\over{2}}(D^{2}-KD)
+\chi  (X,{\cal O}_{X})\cr
&=H^{2}{{n+2}\choose{2}}-{{1}\over{2}}(3H^{2}+KH-2DH){{n+1}\choose{1}}\cr
&+{{1}\over{2}}(H^{2}+KH
-2DH-KD+D^{2})+\chi  (X,{\cal O}_{X}).\cr}
$$
For a rank two vector bundle $E$ the formula reads
$$
\eqalign{
\chi (E+nH)&=H^{2}n^{2}+(c_{1}(E)H-KH)n+{{1}\over{2}}(c_{1}(E)^{2}-Kc_{1}(E))
-c_{2}(E)+2\chi  (X,{\cal O}_{X})\cr
&=2H^{2}{{n+2}\choose{2}}-(3H^{2}+KH-c_{1}(E)H){{n+1}\choose{1}}\cr
&+H^{2}+KH
-c_{1}(E)H-{{1}\over{2}}Kc_{1}(E)+{{1}\over{2}}c_{1}(E)^{2}-c_{2}(E)
+2\chi  (X,{\cal O}_{X}).\cr}
$$

Taking $D=0$ we obtain
$$
\eqalign{
e_{0}(\Oo_{X})&=H^{2}\cr
e_{1}(\Oo_{X})&={{3}\over{2}}H^{2}+{{1}\over{2}}KH\cr
e_{2}(\Oo_{X})&={{1}\over{2}}H^{2}+{{1}\over{2}}KH
+\chi  (X,{\cal O}_{X}),\cr}
$$
and for $D=2K$ 
$$
\eqalign{
e_{0}(\om_{X}^{\otimes 2})&=H^{2}\cr
e_{1}(\om_{X}^{\otimes 2})&={{3}\over{2}}H^{2}-{{3}\over{2}}KH\cr
e_{2}(\om_{X}^{\otimes 2})&={{1}\over{2}}H^{2}-{{3}\over{2}}KH+K^{2}+
\chi  (X,{\cal O}_{X}).\cr}
$$
Finally, taking $E=\Om_{X}$, the cotangent sheaf of $X$, and using the fact
that $c_{1}(\Om_{X})=K$ we deduce 
$$
\eqalign{
e_{0}(\Om_{X})&=2H^{2}\cr
e_{1}(\Om_{X})&=3H^{2}\cr
e_{2}(\Om_{X})&=H^{2}-c_{2}(\Om_{X})+2\chi  (X,{\cal O}_{X}).\cr}
$$

Now the assertion of the remark follows from the Theorem since
$e_{i}(A)=e_{i}(\Oo_{X})$, $e_{i}(\om )=e_{i}(\om_{X} )$ and $e_{i}(\Om )=e_{i}(\Om_{X})
+e_{i}(\Oo_{X})$.\hfill $\Box$
\bigskip

{\bf Theorem 3.4.} Let $k$ be a field, $X\subset {\bf P}^{N}_{k}$ an equidimensional smooth
subscheme of dimension three, $H$ the class of the hyperplane section, $K$ the canonical 
divisor, and $c_{2}$ and $c_{3}$ the second and third Chern class of the cotangent 
bundle of $X$. One has
$$
(H^{3})^{2}\leq 35H^{3}-11H^{2}K-9K^{2}H+c_{2}H-K^{3}-{{1}\over{12}}Kc_{2}
+{{1}\over{2}}c_{3},
$$
and equality holds if and only if $\dim \hbox{\rm Sec}(X)<7$.
\bigskip

{\it Proof.} We use the notation of Theorem 3.1, taking $\ia$ to be generated by 7 general linear forms in
the ideal $I$ of the diagonal. Recall that $\dim {\rm Sec}(X)< 7$ if and only if ${\rm ht}(\ia:I)\geq 8$.

We will apply Theorem 1.9 (a) with $r=s=7$ and $g=4$. Because $X$ is smooth we have
$$\eqalign{
I\omega_{R}/I^{2}\omega_{R} &\build\cong_{7}^{} \omega^{2}\otimes \Omega\cr
I^{2}\omega_{R}/I^{3}\omega_{R} &\build\cong_{7}^{} \omega^{2}\otimes S_{2}\Omega.
}$$

As $X$ is smooth, we can apply
the Riemann-Roch formula as in the proof of Corollary 3.3, to derive the following formulas,
which express the Hilbert coefficients used in Theorem 1.9(a) in terms of the numbers that
appear in our desired formula:

\

$e_0 (A)= H^3$

$e_1 (A)= 2H^3-{{3}\over{2}}KH^2$

$e_2 (A)= {{1}\over{12}} (14H^3+9KH^2+K^2H+c_2 H)$

$e_3 (A)= {{1}\over{24}} (4H^3+6KH^2+2K^2H+2c_2 H+ Kc_2 )$

\ 

$e_0 (\Omega)= 4H^3$

$e_1 (\Omega)= 8H^3+KH^2$

$e_2 (\Omega)= {{1}\over{6}} (28H^3+9KH^2+2K^2H-4c_2 H)$

$e_3 (\Omega)= {{1}\over{12}} (8H^3+6KH^2+4K^2H-8c_2 H+ Kc_2 -6c_3)$

\ 

$e_0 (\om^{\otimes 2})= H^3$

$e_1 (\om^{\otimes 2})= 2H^3-{{3}\over{2}}KH^2$

$e_2 (\om^{\otimes 2})= {{1}\over{12}} (14H^3-27KH^2+13K^2H+c_2 H)$

$e_3 (\om^{\otimes 2})= {{1}\over{24}} (4H^3-18KH^2+26K^2H+2c_2 H-12 K^3-3Kc_2)$

\ 

$e_0 (\om^{\otimes 2}\otimes \Omega)= 4H^3$

$e_1 (\om^{\otimes 2}\otimes \Omega)= 8H^3-7KH^2$

$e_2 (\om^{\otimes 2}\otimes \Omega)= {{1}\over{6}} (28H^3-63KH^2+38K^2H-4c_2 H)$

$e_3 (\om^{\otimes 2}\otimes \Omega)= {{1}\over{12}} (8H^3-42KH^2+76K^2H-8c_2 H-48 K^3+17 Kc_2 -6c_3)$

\ 

$e_0 (\om^{\otimes 2}\otimes S_2\Omega)= 10H^3$

$e_1 (\om^{\otimes 2}\otimes S_2\Omega)= 20H^3-19KH^2$

$e_2 (\om^{\otimes 2}\otimes S_2\Omega)= {{1}\over{6}} (70H^3-171KH^2+119K^2H-25c_2 H)$

$e_3 (\om^{\otimes 2}\otimes S_2\Omega)= {{1}\over{12}} (20H^3-114KH^2+238K^2H-50c_2 H-186 K^3+125 Kc_2 -42c_3)$
\hfill$\Box$
\bigskip

{\bf Corollary 3.5.} {\it Let $k$ be a perfect field, $X\subset {\bf P}^{N}_{k}$ an 
equidimensional smooth subscheme of dimension three, $A$ 
its homogeneous coordinate ring, $\om :=\om_{A}$ the canonical module, and 
$\Om :=\Om_{A/k}$ the module of differentials. One has
$$
e_{0}(A)^{2}+391e_{0}(A)-246e_{1}(A)+66e_{2}(A)+50e_{3}(A)\geq 18e_{2}(\om \otimes_{A}\om )
-2e_{3}(\Om )-2e_{3}(\om \otimes_{A}\om )
$$
and equality holds if and only if $\dim\hbox{\rm Sec}(X)<7$.}

{\it Proof.} Use the formulas given in Theorem 3.4 to express all the necessary Hilbert coefficients
in terms of $e_{0}, \dots, e_{3}$. 
\hfill$\Box$
\bigskip

{\bf Remark 3.6.} If $H$ is the class of the hyperplane section, $K$ the canonical divisor,
$D=H^3$ the degree of $X$, the above inequality is equivalent to:
$$
D^{2}\geq 7(5D+3KH^{2}+K^{2}H-c_{2}(\Om_{X})H)-2c_{2}(\Om_{X})K+K^{3}+c_{3}(\Om_{X}).
$$
In other words, if $C$ and $S$ are respectively a curve and a surface obtained
by taking general linear sections of $X$, the formula reads
$$
D^{2}\geq 7(5D+3\chi_{C}+12\chi ({\cal O}_{S})-2\chi_{S})-48\chi ({\cal O}_{X})+K_{X}^{3}+\chi_{X}.
$$

\bigskip\bigskip

{\bigrm References}\bigskip

[1]{} R. Ap{\'e}ry, Sur les courbes de premi{\`e}re esp{\`e}ce de
l'espace {\`a} trois dimensions, {\sl C. R. Acad. Sci. Paris}, t. 220,
S{\'e}r. I, p. 271--272, 1945. 
\medskip

[2]{} M. Artin and M. Nagata, Residual intersections in Cohen-Macaulay
rings, {\sl J. Math. Kyoto Univ.} {\bf 12} (1972), 307--323.
\medskip

[3]{} L. Avramov and J. Herzog, The Koszul algebra of a codimension 2
embedding, {\sl Math. Z.} {\bf 175} (1980), 249--280.
\medskip

4] R. Buchweitz, Contributions \`a la th\'eorie des singularit\'es, These d'Etat, Universit\'e Paris 7, 1981.
Available at https://tspace.library.utoronto.ca/handle/1807/16684.
\medskip

[5] H. Cartan and S. Eilenberg, Homological Algebra, Princeton University Press, Princeton, NJ.
\medskip

[6]{} M. Chardin, D. Eisenbud, and B. Ulrich {\sl Hilbert functions and
residually $S_{2}$ ideals}, Compositio {\bf 125} (2001), 193--219. 
\medskip

[7] C. Cumming, Residual intersections in Cohen-Macaulay rings. J. Algebra 308 (2007), no. 1, 91--106.
\medskip

[7a] M. Dale,
Severi's theorem on the Veronese-surface. 
J. London Math. Soc. 32 (1985) 419--425.
\medskip

[8] D. Eisenbud and E. G. Evans, Generating modules efficiently: theorems from algebraic K-theory. J. Algebra (1973), 278--305. 
\medskip

[9]{} F. Gaeta, Quelques progr{\`e}s r{\'e}cents dans la classification
des vari{\'e}t{\'e}s alg{\'e}briques d'un espace projectif, {\sl
Deuxi{\`e}me Colloque de G{\'e}ometrie Alg{\'e}brique}, Li{\`e}ge, 1952. 
\medskip

[10] S. H. Hassanzadeh, Cohen-Macaulay residual intersections and their Castelnuovo-Mumford regularity. Trans. Amer. Math. Soc. 364 (2012) 6371--6394. 
\medskip

[11]  J. Herzog, Komplexe, Aufl\"osungen und Dualit\"at in der lokalen
Algebra. Habilitationsschrift (1974). 
\medskip

[12]{} J. Herzog, A. Simis, and W.V. Vasconcelos, Koszul homology and
blowing-up rings, in {\sl Commutative Algebra}, eds. S. Greco and
G. Valla, Lecture Notes in Pure and Appl. Math. {\bf 84}, Marcel Dekker,
New York, 1983, 79--169. 
\medskip

[13]{} J. Herzog, W.V. Vasconcelos, and R. Villarreal, Ideals with
sliding depth, {\sl Nagoya Math. J.} {\bf 99} (1985), 159--172.
\medskip

[14]{} C. Huneke, Linkage and Koszul homology of ideals, {\sl Amer. J.
Math.} {\bf 104} (1982), 1043--1062.
\medskip

[15]{} C. Huneke, Strongly Cohen-Macaulay schemes and residual
intersections, {\sl Trans. Amer. Math. Soc.} {\bf 277} (1983), 739--763.
\medskip

[16]{} C. Huneke and H. Martin, Residual Intersection and the Number of
Equations Defining Projective Varieties, Comm. Algebra
{\bf 23} (1995), no. 6, 2345--2376.
\medskip

[17] A. Simis and B. Ulrich, On the ideal of an embedded join. J. Algebra 226 (2000), no. 1, 1Ð14.

[18] J. St{\"u}ckrad, On quasi-complete intersections, {\sl Arch. Math.}
{\bf 58} (1992), 529--538.\medskip 

[19]{} B. Ulrich, Remarks on residual intersections, in {\sl Free
Resolutions in Commutative Algebra and Algebraic Geometry, Sundance 1990},
eds. D. Eisenbud and C. Huneke, Res. Notes in Math. {\bf 2}, Jones
and Bartlett Publishers, Boston-London, 1992, 133--138.
\medskip

[20]{} B. Ulrich, Artin-Nagata properties and reductions of ideals,
{\sl Contemp. Math.} {\bf 159} (1994), 373--400.
\medskip

[21]{} J. Watanabe, A note on Gorenstein rings of embedding
codimension three, {\sl Nagoya Math. J.} {\bf 50} (1973), 227--232.
\medskip

\bigskip\bigskip

\noindent Marc Chardin, Institut de Math{\'e}matiques, 
CNRS \&\ Universit{\'e} Pierre et Marie Curie

\qquad chardin@math.jussieu.fr

\bigskip

\noindent David Eisenbud, 
Department of Mathematics, University of California, Berkeley 

\qquad de@msri.org

\bigskip

\noindent Bernd Ulrich, Department of Mathematics, Purdue
University, 

\qquad  ulrich@math.purdue.edu 

\end